# Safety stock calculation with heteroscedastic, autocorrelated and heavy-tailed forecast errors, an adaptive parametric method.

Keywords: Safety stock; Forecast error; Autocorrelation; Heteroscedasticity; Heavy tails; Fill rate.

**Authors:** Luis Fernández-Palacios [a,*], Manuel Ceballos [b], Yolanda Muñoz-Ocaña [c]

[a] *Loyola University Andalusia, Department of Business Management, Campus Sevilla: Avda. de las Universidades, 2, 41704 - Dos Hermanas, Sevilla, Spain. E-mail address: lfpalacios@uloyola.es*

[b] *Loyola University Andalusia, Department of Engineering, Campus Sevilla: Avda. de las Universidades, 2, 41704 - Dos Hermanas, Sevilla, Spain. E-mail address: mceballos@uloyola.es*

[c] *Loyola University Andalusia, Department of Business Management, Campus Córdoba: Calle Escritor Aguayo, 4, 14004 Córdoba, Spain. E-mail address: yolandam@uloyola.es*

[*] *Corresponding author. Calle Estrellita Castro, 5, 41002 – Sevilla, Spain. E-mail address: lfpalacios@uloyola.es*

***Data statement:*** *The data used in this study are subject to commercial confidentiality restrictions and are therefore not publicly available. The dataset contains sensitive business information from a private company. The key replicable components of the study are the mathematical derivations, assumptions, and analytical procedures, which are fully described in the manuscript.*

# Abstract

This study develops and validates an adaptive parametric method for sizing safety stock when forecast errors are autocorrelated, heteroscedastic, asymmetric, and heavy-tailed. Using four years of weekly forecasts from a fast-moving consumer goods company, the analysis documents substantial positive autocorrelation, widespread departures from normality, systematic forecast optimism, and heterogeneous relationships between forecast volume and error variance.

The contribution is a modular, adaptive, and sequential safety-stock architecture that integrates corrections within a single SKU-level protection-point calculation. The architecture first constructs cumulative forecast errors over the operational action period, aligning uncertainty measurement with replenishment decisions while preserving temporal dependence. It then selects an SKU-specific parametric error distribution and determines a baseline protection point through an expected-shortage condition consistent with the target beta fill rate. From this baseline, a forecast-scaled module adjusts distributional dispersion when error variance depends on forecast volume, and a Cycle-End Exposure safeguard subsequently limits within-cycle shortage risk when large replenishment lots create additional exposure.

Performance is evaluated through discrete-event simulation against the company's current method. The proposed method reduces average inventory by 9.9% and maximum inventory by 19.5%, while showing no statistically significant deterioration in beta fill rate or average stockout duration. The proportion of weeks without stockouts declines slightly, without evidence of more persistent or operationally severe shortages. The heteroscedasticity module provides selective benefits, with larger service improvements as variance dependence intensifies.

The resulting framework is interpretable, adapts to SKU-specific statistical profiles, preserves traceability between distributional features and safety-stock outcomes, and can be integrated into existing ERP environments.

# 1. Introduction

This article presents a real case study of a Fast-Moving Consumer Goods (FMCG) company and proposes an alternative to the method currently used for safety stock (SS) sizing. In a competitive environment that demands continuous cost reduction, the company aims to reduce inventory levels without compromising the target service level. Previous internal analyses had estimated that supply lead-time variability accounted for less than 4% of total variance, indicating that it was a comparatively minor source of uncertainty in the company's operating context. Accordingly, supply uncertainty is treated as outside the scope of the present study, which focuses on forecast-error uncertainty as the main driver of safety-stock requirements. This focus is also consistent with the emphasis observed in the inventory literature: Barros et al. (2021) found that only 10% of the reviewed studies addressed supply uncertainty exclusively, whereas 30% considered it jointly with demand uncertainty.

The main advantage of using forecasts in production planning lies in improving the temporal positioning of supply decisions. Unfortunately, "All that we can say for certain about a forecast of demand is that it will be in error" (Silver et al., 1998, p. 109). Forecast-error uncertainty is a central driver of inventory-related costs: when forecast errors are not properly modelled, they propagate upstream through production, replenishment, and capacity-planning decisions (Saoud et al., 2025). The pandemic altered purchasing patterns so profoundly that naïve forecast errors increased significantly (E2open, 2024), prompting firms to adopt inventory practices based on larger buffer stocks, thereby tying up substantial working capital (Song et al., 2025). This gives rise to what Wang et al. (2024) call "algorithm aversion", whereby organisations or individuals reject optimal analytical decision support in favour of informal and subjective decisions, even though algorithmic decisions may significantly outperform human judgment.

Each model studied in the scientific literature has limitations determined by the assumptions it makes beforehand about reality. Graves (1999) pointed out that one of the biggest constraints on the continued development of inventory theory is how to incorporate more realistic assumptions about product demand into new models. In their reference works, Chopra & Meindl (2016, Chapter 12) and Silver et al. (1998, Chapter 3) highlight the limitations of theoretical methods and they require empirical validation to assess their applicability under operational conditions. Unfortunately, real and representative data is not always available for researchers, of the more than 300 articles analyzed, less than 15% are supported by field data.

The company's current method for calculating SS follows the classical inventory-theory formulation based on forecast-error uncertainty. Because this study focuses exclusively on forecast-error variability, lead-time uncertainty is excluded from the expression. The data used are weekly, as is the planning process.

$$SSweek = Z\sigma_e\sqrt{L}$$

Where

SSweek: Weekly safety stock

Z: Service-level factor

$\sigma_e$: Weekly demand forecast error standard deviation

L: Average lead time

The company expresses this estimate as a fixed number of days of forecast coverage. For each planning period, SS in units is calculated by applying these coverage days to the corresponding forecast rate and therefore varies proportionally with the forecast.

Chatfield (2000) strongly criticized this expression on the grounds that it has no theoretical justification and can yield very poor results. The algorithm is based on the additive property of variance for independent random variables. For this relationship to be valid, errors in each period must:

I. Be independent (iid). If the forecast error exhibits temporal dependence, there is autocorrelation between periods, then the variance of cumulative consumption over a lead time does not simply scale linearly with L.
II. Have a variance that is a representative measure of dispersion. Under heavy tails or strong asymmetry, the standard deviation may be dominated by extreme observations and fail to reflect the operational risk relevant for SS.

III. Be centred around zero, i.e., forecast errors must not contain systematic bias. The classical formulation assumes errors centred around zero; otherwise, bias and random variability are mixed, potentially distorting safety-stock estimates.

IV. Have the same variance in each period. The company relaxes this assumption by scaling SS proportionally with the weekly forecast, implicitly assuming the same proportional relationship for forecast-error dispersion.

Additionally, the classical formulation targets the probability of completing a replenishment cycle without a stockout. However, this is not the metric by which customers evaluate our performance. In the company under study, the relevant service indicator is the beta fill rate (BFR), defined as the ratio between units served and units requested.

Therefore, the current procedure relies on a set of assumptions that, although embedded in the design of the algorithm, are not verified in practice. As a result, the algorithm is repeatedly adjusted through informal and subjective decisions. The objective of this study is to:

I. Develop a method capable of overcoming the current limitations, in order to avoid subsequent manual adjustments based on subjectivity and, consequently, reduce so-called “algorithm aversion”.

II. Establish a clear relationship between the final calculation and its determining variables. This last aspect is especially relevant in an operational environment, since, as Silver et al. (1998, Chapter 3) note, lack of understanding by management decision-makers can become an important obstacle to the implementation of new models.

III. Adapt the method to the statistical profile of each SKU, avoiding the application of a single uniform solution across heterogeneous cases.

IV. Ensure compatibility with the company’s current ERP while maintaining sufficient stability to avoid continuous maintenance and recalculation.

Existing approaches can be broadly classified into normal approximations, correction-factor methods, non-parametric techniques, and explicit modelling of the forecast-error distribution. This study builds on the latter approach by developing an adaptive representation of forecast uncertainty that accounts for its statistical and operational characteristics and translates them into the resulting protection level. To the best of our knowledge, no previous study has proposed a unified parametric method that jointly addresses the four identified limitations, adapts to the statistical profile of each SKU, and preserves operational efficiency.

# 2. Literature review: Study of Limitations

## 2.1 Autocorrelation of forecast errors

Significant levels of autocorrelation have been identified in demand series obtained from real business data (Erkip & Hausman, 1994; Lee et al., 2000). This dependence between consecutive observations has important operational consequences, as stockouts tend to become both larger and more frequent when demand observations are correlated (Ray, 1980; Zinn et al., 1992). Moreover, lead-time forecasts can be interpreted as a form of overlapping temporal aggregation, which may further reinforce the dependence between consecutive forecast observations and forecast errors (Boylan & Babai, 2016).

Under these conditions, the conventional square-root-of-lead-time relationship no longer provides an adequate representation of uncertainty. When demand is autocorrelated, SS may increase faster than the square root of lead time (Graves, 1999). The resulting inventory and service costs have been widely documented in the literature (Díaz, 2007; Syntetos & Boylan, 2008; Díaz & Behr, 2014; Boute et al., 2014; Prak et al., 2017).

A theoretically exact estimation of this effect would require the covariance structure between consecutive weekly forecast errors to be known for each SKU. In practice, however, estimating these covariance factors reliably at the individual-item level is complex. To address this difficulty, a cumulative smoothed-MSE procedure has been proposed for estimating lead-time demand variability (Syntetos & Boylan, 2006). A widely used alternative is to calculate the empirical variance of cumulative forecast errors, thereby incorporating the dependence structure without estimating each covariance term separately (Eppen & Martin, 1988; Rostami-Tabar et al., 2013, 2014; Lee, 2014; Trapero et al., 2019a; Derbel et al., 2022). This empirical approximation performs well in both simulated and real-data settings (Saoud et al., 2022), and subsequent research has extended and confirmed these findings under stochastic lead times (Babai et al., 2022).

However, lead-time forecast errors are typically defined over the cumulative lead-time horizon without considering when the forecast was originally produced. This distinction is operationally relevant because the inventory physically available in period (t+i) results from replenishment decisions made approximately (i) periods earlier, where (i) corresponds to the replenishment lead time (L). Because forecast accuracy generally improves as the forecast horizon shortens, ignoring this operational lag can lead to underestimation of the forecast-error variance relevant to the replenishment decision and, consequently, to underestimated safety-stock requirements. Our work addresses this limitation by aligning forecast-error measurement with the timing of the replenishment decision, using the information that was genuinely available when the order was placed.

## 2.2 Heavy-tailed forecast-error distributions

Empirical evidence indicates that the normality assumption is frequently violated in real inventory applications. Both lead-time demand distributions (Tyworth & O'Neill, 1997; Vernimmen et al., 2008) and forecast-error distributions (Coleman & Swanson, 2007) may exhibit non-normal and asymmetric behaviour. Forecast errors are typically concentrated around relatively small deviations, but a limited number of substantially larger errors generate long-tailed distributions (Spiliotis et al., 2019). This pattern is reflected in the fact that extreme forecast errors were almost four times greater for tail items than for high-rotation items (E2open, 2019).

The normal approximation may nevertheless perform reasonably well under certain conditions, particularly when the coefficient of variation is below 0.3 or when lead times are sufficiently long (Syntetos & Boylan, 2008). These conditions are not fully applicable to our setting, where lead times are limited to four weeks. At the individual-item level, moreover, assuming normally distributed forecast errors tends to underestimate the actual degree of uncertainty (Lee & Scholtes, 2014; Barrow & Kourentzes, 2016).

Classical inventory methods may perform poorly when demand or forecast errors are non-normal because imposing normality distorts inventory decisions and costs (Mentzer & Krishnan, 1985; Eppen & Martin, 1988; Lau & Lau, 2003; Das et al., 2014; Silver et al., 2016). Parametric alternatives select SKU-specific distributions and use their quantiles to determine inventory levels, as demonstrated for intermittent spare-parts demand (Turrini & Meissner, 2019). Other work evaluates forecasting methods through empirical accuracy distributions but does not translate tail risk into safety-stock optimisation (Spiliotis et al., 2019).

Beyond these parametric and empirical approaches, a substantial body of research has addressed non-normal demand by relaxing distributional assumptions. Distribution-free inventory control has been proposed for highly uncertain demand (Gallego et al., 2007), while other contributions account for heterogeneous or multimodal forecast-error distributions requiring specific inventory corrections (Manary et al., 2009). Scenario-based approaches estimate service levels from possible lead-time demand outcomes (Ruiz-Torres & Mahmoodi, 2010), and linear-programming formulations have been proposed when the lead-time demand distribution is only partially known (Janssens & Ramaekers, 2011).

The persistent lack of normality has also encouraged the development of non-parametric and tail-sensitive methods. The Expected Value Reorder-point model (EVR) derives reorder-point decisions from possible lead-time demand outcomes rather than from a normal approximation (Ruiz-Torres & Mahmoodi, 2010). Extreme value theory explicitly models tail behaviour associated with unusually large observations (Biçer, 2015), while kernel density estimation (KDE) flexibly represents the empirical distribution without imposing a predefined parametric form (Trapero et al., 2019a). Conditional Extreme Value Theory (CEVT) extends this focus to conditional tail behaviour (Derbel et al., 2022), and bootstrap-based approaches have been applied to safety-stock estimation under skewed, multimodal, or highly variable demand (Saldanha et al., 2023).

Overall, four broad responses to non-normality can be identified in the literature: assuming normality, applying correction factors, using non-parametric methods, and directly modelling the shape of the observed distribution. Our research adopts the fourth approach. Its main conceptual advantage is that the adjustment is derived from the same model used to generate the inventory parameters, thereby preserving the traceability of the underlying cause-and-effect relationships and supporting managerial decision-making. Its main limitation is the risk of distributional misspecification. In the present setting, this risk is mitigated by the available historical depth, as discussed in Section 3. To the best of our knowledge, the application of this approach to forecast-error distributions in fast-moving consumer goods remains largely unexplored.

## 2.3 Bias and asymmetric forecast-error distributions

When forecast-error distributions are asymmetric, probability mass may be incorrectly allocated across the two sides of the distribution, leading to substantial errors in safety-stock estimation (Clark, 1957; Krupp, 1997). In the presence of optimistic bias, this misallocation can artificially inflate the calculated SS. A slight positive bias may sometimes help achieve the target service level, although at the cost of higher inventory (Syntetos & Boylan, 2008). Its effect, however, is neither linear nor symmetric. When positive bias remains in the error sample used to estimate SS, inventory may be inflated twice: first through systematically overstated forecasts and again by treating the same bias as random uncertainty (Syntetos & Boylan, 2008).

Bias is also closely linked to the way forecasts are produced and revised in practice. Judgmental adjustments are common, with approximately 45% of planners reporting manual modifications to statistical forecasts (Sanders & Manrodt, 1994). These interventions may introduce behavioural bias that is difficult to observe directly as a causal variable. At the same time, human revisions are not necessarily detrimental, as they may improve forecast accuracy under certain conditions (Mathews & Diamantopoulos, 1992; Fildes et al., 2009).

The literature has addressed this problem through several complementary mechanisms. Planned lead times and SSs have been used to mitigate the combined effects of forecast bias and demand uncertainty in MRP-controlled batch-production systems (Enns, 2001). Forecasting practice has also emphasised that bias should be addressed before forecasts are incorporated into decision models (Crum & Palmatier, 2003). Specific adjustment procedures have been proposed for biased sales forecasts (Manary & Willems, 2008), while industrial applications have shown how forecast-error bias and heterogeneity can be corrected in practice (Manary et al., 2009). More recently, the classical Martingale model has been extended by explicitly introducing a bias parameter in a simulation setting (Seiringer et al., 2022).

As discussed above, ignoring bias in the forecast-error sample creates risks of both excess inventory and stockouts. Correcting it, either in the forecast itself or in the resulting forecast error, is nevertheless difficult because the structural component of bias must be distinguished from the portion that lies beyond the model's predictive capacity. Our proposed method addresses this limitation by fitting the complete forecast-error distribution without imposing a zero-mean restriction and, when supported by the data, allowing for asymmetric distributional shapes. The resulting location and asymmetry parameters are then incorporated directly into the expected-shortage calculation, making their effect on the protection point explicit and traceable.

## 2.4 Conditional heteroscedasticity in forecast errors

In practice, demand and forecast errors are influenced by exogenous factors such as prices, weather conditions, promotions, and calendar effects, which limits the validity of classical assumptions based on stationary and homogeneous uncertainty (Beutel & Minner, 2012). Safety-stock requirements may increase substantially when demand is not generated by a stationary independent and identically distributed process (Erkip et al., 1990). Similarly, when forecast-error variance changes over time, a constant safety-stock level may overprotect low-risk periods while underprotecting high-risk periods (Campbell, 1995). Heteroscedasticity appears widespread, with previous studies detecting it in more than 80% of the industries analysed (Kanet et al., 2010) and in 85% of monthly EUROSTAT series, thereby supporting period-specific SS (Stößlein et al., 2014).

Several approaches have been developed to account for this variation in uncertainty. Time-Based Mean Absolute Deviation provides a time-sensitive measure of forecast error and is similar to the method used by the company analysed in our study (Krupp, 1997). Temporal heteroscedasticity has also been identified in industrial demand series, although its implications have historically received limited attention in the inventory-management literature (Zhang & Kline, 2007). Forecast-conditional kernel smoothing has been used to represent heterogeneous and multimodal forecast errors associated with sales peaks and promotions (Manary et al., 2009). Causal variables such as price, weather, and weekday have been incorporated into safety-inventory planning through OLS regression and linear-programming approaches (Beutel & Minner, 2012).

More recent contributions model heteroscedasticity as a non-constant conditional variance. GARCH(1,1) models have been applied to represent time-varying forecast-error volatility (Trapero et al., 2019a, 2019b), while ARIMA-GARCH has been used to model monthly cement demand in India, although without reporting its impact on safety-stock reduction or service levels (Jaipuria & Mahapatra, 2021). Other approaches combine GARCH with Filtered Historical Simulation and Conditional Extreme Value Theory to address heteroscedasticity, non-normality, and extreme demand jointly (Derbel et al., 2022).

The literature identified therefore concentrates primarily on temporal heteroscedasticity, where variance changes across periods. The sizing of SS remains less studied when conditional variance depends explicitly on observable covariates, and only one identified study considers forecasted sales intensity as the covariate driving that variance, without jointly incorporating autocorrelation, non-normality, and bias. Our study addresses these gaps through a practical method that links conditional variance to forecasted sales intensity and simultaneously accounts for the other distributional limitations within a single framework that remains interpretable for managerial decision-making.

## 2.5 Integration gap and architectural contribution

The literature provides several statistical and inventory-theoretic elements for addressing individual violations of the classical safety-stock assumptions; however, these solutions remain fragmented. To the best of our knowledge, no previous study combines temporal dependence, non-normality and bias, forecast-dependent heteroscedasticity, and lot-size-related exposure within a single parametric framework validated with real operational data. This study addresses that integration gap by developing a unified SKU-level method that links statistical diagnosis to a BFR-based operational protection-point calculation.

# 3. Empirical evidence of violated assumptions

Six years of forecast records were available, with weekly forecasts for the subsequent 58 weeks for each SKU. For Hyndman and Kostenko (2007), calendar years are relevant in a seasonal business; the historical window used should cover multiples of calendar years to include all seasonal events. We exclude the first year, 2020, from the analysis because it coincides with the onset of the COVID-19 pandemic, when market dynamics were severely disrupted. Chatfield (2000, Chapter 6) warns us that historical windows including structural breaks can bias model parameters and it is often preferable to restrict the historical window to post-change data. According to the E2open (2024) Benchmark, at the peak in the Spring of 2020, almost half of all volume was exposed to extreme error. We reserve the final year, 2025, for out-of-sample validation through discrete-event simulation, leaving four complete calendar years, 2021–2024, for model development and parameter estimation. This calibration window was previously found to provide a representative basis for estimating forecast-error uncertainty at high service levels on the same dataset (Fernández-Palacios et al., 2026).

Temporal dependence is assessed using the Durbin–Watson statistic and the first-order autocorrelation coefficient, $\rho_1$. Values of $DW < 1.5$ are classified as evidence of strong positive persistence, while $DW > 2.5$ indicates strong negative alternation. Strong positive autocorrelation is present in 54% of the SKUs, whereas strong negative autocorrelation appears in only 4%. Within the group with $DW < 1.5$, the mean first-order autocorrelation is $\rho_1 = 0.39$, with values reaching 0.8, comparable to those reported by Erkip and Hausman (1994) and Lee et al. (2000). These results indicate that forecast errors are not independent for a substantial proportion of the portfolio.

Throughout the empirical analysis, forecast error is defined as observed sales minus forecast. Therefore, negative errors indicate forecast optimism, whereas positive errors indicate demand exceeding the forecast. The normal approximation is first evaluated using the coefficient-of-variation threshold of 0.3 proposed by Syntetos and Boylan (2008). In the sample, 99.0% of SKUs exceed this threshold, suggesting that the normal approximation may be unsuitable for most of the portfolio. The Shapiro–Wilk test rejects normality for 40.1% of SKUs. Extreme skewness and kurtosis affect 8.0% of the sample, but asymmetry remains operationally relevant. Among SKUs with extreme skewness, 76.4% exhibit pronounced forecast optimism. Even among SKUs with less extreme asymmetry, forecast optimism remains the dominant pattern, representing 61.6% of cases. Therefore, forecast errors frequently depart from a centred normal distribution. This is consistent with the empirical evidence reported by Coleman and Swanson (2007), Lee and Scholtes (2014), and Barrow and Kourentzes (2016), who show that real forecast-error distributions are often non-normal, asymmetric, and heavy-tailed.

Heteroscedasticity may arise either across forecast horizons or conditionally on external variables. In this study, the first source is not central because forecasts are available for horizons longer than the replenishment lead time, allowing uncertainty to be measured directly over the operational action period relevant to the replenishment decision. We therefore focus on forecast volume as the conditioning variable because it is directly available within the planning process at the required SKU-week granularity, whereas other potentially relevant drivers may be unavailable, fragmented, or not known with sufficient accuracy at the time the replenishment decision is made. Higher forecast volumes may also reflect greater commercial intensity and therefore be associated with different forecast-error variance. Conditional heteroscedasticity is assessed using the Breusch–Pagan test, with $p < 0.05$ as the rejection criterion. Evidence of heteroscedasticity is found in 46% of the portfolio. The intensity of the relationship between forecast volume and error variance is also highly heterogeneous. Among the SKUs for which the null hypothesis is rejected, the Breusch–Pagan Lagrange Multiplier statistic ranges approximately from 4 to 55. Thus, forecast-error variance does not increase proportionally with forecast volume across the entire portfolio, nor does the relationship have a uniform magnitude among affected SKUs. This result is consistent with studies such as Kanet et al. (2010), Stößlein et al. (2014), and Trapero et al. (2019a), which report substantial heterogeneity in time-varying demand and forecast-error variance.

The empirical evidence does not support the assumptions of independent, normally distributed, unbiased forecast errors, nor a uniformly proportional relationship between forecast volume and forecast-error variance across the portfolio. These findings justify a safety-stock method that adapts to the statistical profile of each SKU rather than applying a uniform correction to the entire portfolio.

# 4. Proposed mathematical method

The proposed mathematical method is structured as a sequence of interdependent modules. Cumulative forecast error in the action period (CFEAP) first defines the operational forecast-error variable by aligning cumulative errors with the information available when the replenishment decision was made. Adaptive Parametric Fill-Rate Method (APFR) then identifies the SKU-specific parametric distribution of that error and translates it into a protection point consistent with the target BFR. For SKUs with forecast-dependent error variance, Forecast-Scaled Error Distribution Method (FSEDM) conditionally rescales the fitted distribution according to the forecast level, with the magnitude of the adjustment determined by the estimated intensity of heteroscedasticity. Finally, the Cycle-End Exposure rule constrains solutions that may satisfy the average fill-rate objective but provide insufficient protection near the end of long replenishment cycles. The output of each module therefore becomes an input or constraint for the next, producing a single coherent SKU-level safety-stock decision.

## 4.1 Autocorrelation: CFEAP

As discussed in Section 2.1, the use of cumulative forecast errors is established in the literature. The remaining limitation is that the error horizon is not always linked to the forecast origin that informed the replenishment decision. CFEAP addresses this limitation by aggregating the forecast errors corresponding to the information set available when the production order was placed. Its contribution is therefore operational alignment rather than cumulative aggregation per se, making the measure more appropriate for MRP and production-planning environments. For each week t, we define the cumulative forecast error as the sum of the differences between real demand and forecast from t+L to t+2L−1. This represents the relevant error for decisions that are executed after the lead time. The empirical variance of this series is used as an estimator of risk over the control horizon.

$$E_t = \sum_{i=L}^{2L-1} (y_{t+i} - \hat{y}_{t+i|t})$$

Where:

3. L = Lead time
4. t+i: real demand in week t+i
5. $\hat{y}_{t+i|t}$: forecast made in week t for week t+i
6. $E_t$: total error over the L future periods where you can still take action

The resulting CFEAP series consists of overlapping windows, each associated with a distinct replenishment decision and forecast origin. To illustrate this structure, consider a lead time of $(L = 4)$ weeks. The cumulative forecast error associated with the decision made in week $(t)$ is:

$$C_t = e_{t,t+4} + e_{t,t+5} + e_{t,t+6} + e_{t,t+7}$$

For the following weekly decision, the cumulative forecast error becomes:

$$C_{t+1} = e_{t+1,t+5} + e_{t+1,t+6} + e_{t+1,t+7} + e_{t+1,t+8}$$

Where $e_{t,s}$ denotes the difference between realized demand in week (s) and the forecast for that week issued at forecast origin (t).

Although $C_t$ and $C_{t+1}$ represent two distinct operational decisions, their evaluation horizons share three weeks of realized demand $t+5, t+6, and\ t+7$. However, the corresponding forecast errors are not identical, since they are calculated using forecasts issued at different origins. This overlap is not a limitation specific to CFEAP, but an inherent consequence of representing weekly replenishment decisions over overlapping future horizons.

By working directly with the empirical distribution of cumulative forecast errors, rather than decomposing risk into individual weekly variances, CFEAP captures the observed distributional characteristics of the cumulative errors, including asymmetry and tail behavior. Nevertheless, overlapping evaluation horizons may generate temporal dependence between consecutive observations. The extent of this dependence depends on the underlying demand process and forecast revisions, and cannot be inferred from the degree of overlap alone.

Consequently, CFEAP observations should not be assumed to be independent and identically distributed (i.i.d.). Temporal dependence may reduce effective sample size and affect uncertainty estimation, particularly for tail quantiles, without invalidating the empirical distribution for replenishment-risk assessment. Therefore, the observed dependence structure should be considered rather than assuming an effective sample size of $n/L$. A subsequent empirical study (Fernández-Palacios et al., 2026) further characterized the temporal dependence of CFEAP observations and employed a Moving Block Bootstrap to account for its effects on estimation uncertainty.

## 4.2 Non-normality and Bias: APFR

APFR is the integrated module that replaces the normality assumption. Its contribution lies in coupling parametric distribution fitting with expected-shortage optimisation within a common SKU-specific decision rule. Statistical fit alone is not considered sufficient for model selection. Candidate distributions must also generate protection points supported by the observed right-tail information and consistent with operational service performance. The resulting fill-rate solution is subsequently constrained by the Cycle-End Exposure safeguard, which limits expected shortage using an operational weekly-demand threshold.

### 4.2.1 Probability density function identification through BDFT

A parametric model selection procedure is required to identify, for each SKU, the probability distribution that best characterizes its probability density function (PDF), instead of imposing a single distributional assumption across all cases. For this purpose, we apply a Best-Fitting Distribution Test, a statistical procedure designed to determine which theoretical probability distribution best fits a given set of observed data. The procedure considers a set of candidate distributions and evaluates their relative fit using the AIC.

A fitted distribution was considered admissible when the protection-point equation for the target BFR produced a finite solution within a predefined data-supported range. Distributions requiring solutions far beyond the observed error range were excluded because the resulting protection points depended primarily on unsupported tail extrapolation rather than on the empirical information available for the SKU. A second admissibility condition required the fitted distribution to provide adequate protection in the operationally relevant right-tail region.

For several fitted families, particularly Weibull minimum, Cauchy, and Lognormal, the protection point required to achieve the 98% BFR was located far beyond the range supported by the observed forecast errors. Although a mathematical solution could sometimes be obtained by extending the numerical search interval, the resulting SS was driven by extreme tail extrapolation and was therefore considered insufficiently supported by the available data. These fits were classified as inadmissible for operational use. Gamma fits were treated separately. Although their protection-point solutions were numerically feasible, they systematically underestimated the operationally relevant right-tail risk. In the validation analysis, the resulting protection levels were insufficient to achieve the target service performance.

Among the admissible candidates, Fitted Beta, Normal, and Logistic provided the best representation of forecast-error uncertainty for 48.4%, 45.1%, and 6.5% of the SKUs, respectively (Figure 1). Beta distributions were fitted using shape parameters, together with location and scale parameters. The flexibility of the Beta distribution allows asymmetric error profiles to be represented without imposing symmetric protection. It also does not assume heavy tails, since it has bounded support.

The Logistic distribution is symmetric like the Normal distribution but has heavier tails, providing greater right-tail protection for SKUs with extreme deviations for which the Normal distribution underestimates tail risk. Unlike Beta, Logistic distribution remains symmetric and therefore cannot represent distributional skewness. It is also numerically stable and always yields a solution for the fill rate. In the Normal case, our approach still improves on the basic formula by capturing bias rather than assuming a zero-mean error distribution.

**Figure 1. Effect of Distributional Shape on the 98% Protection Point**

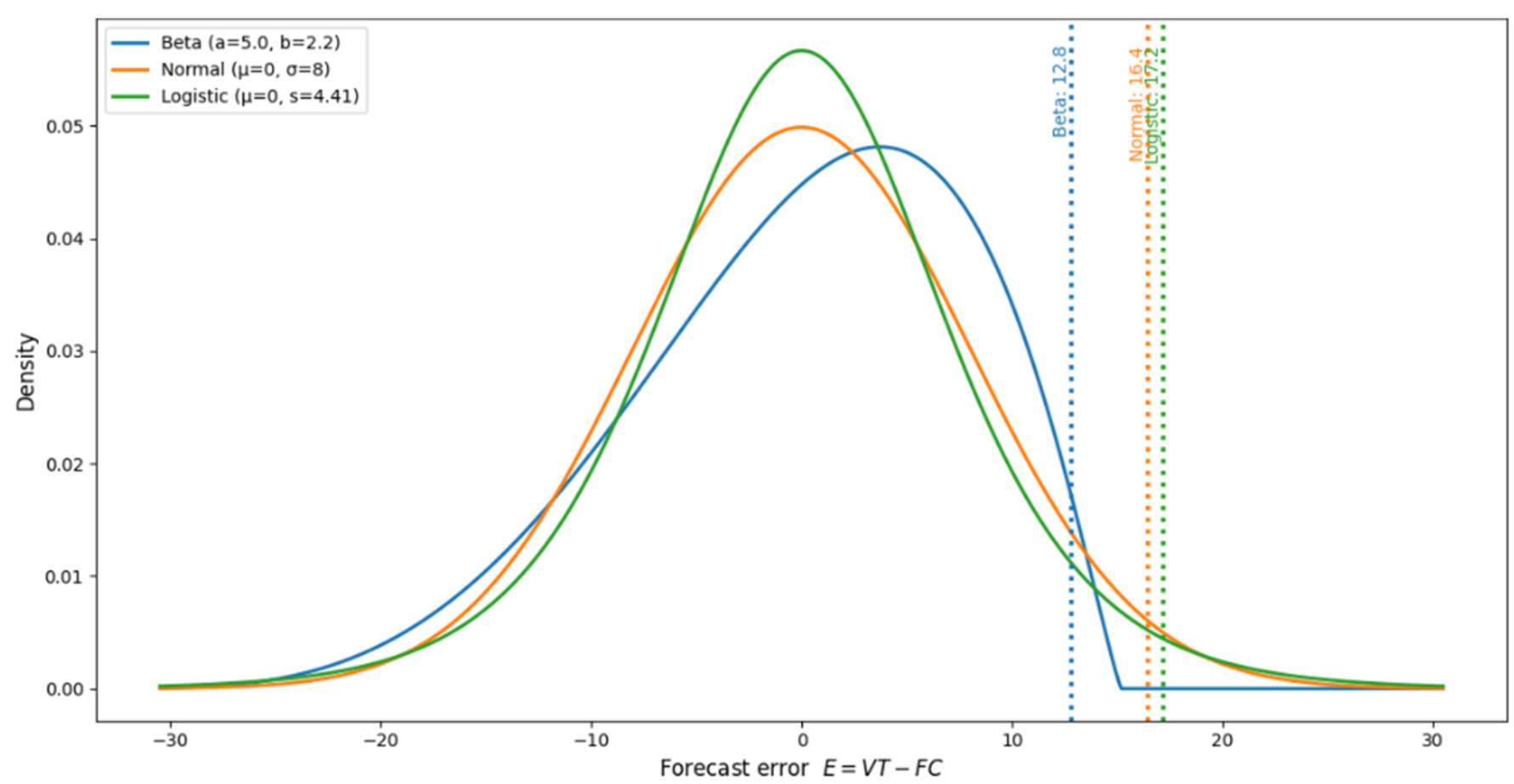

Once these three distributions were selected, out of the 218 SKUs sold in the year in which we ran the DES, we were only able to reach an admissible complete solution for 153. Although only 153 of the 218 SKUs satisfied the inclusion criteria, these SKUs accounted for 96% of production volume. The three main reasons for excluding the remaining SKUs from the final analysis were:

1. **Lack of historical records:** We could not determine a distribution via BDFT for 16 SKUs because of insufficient records for newly introduced products.

2. **Strongly negatively biased forecasts:** When most of the probability mass lies below zero, reflecting systematic forecast optimism, the expected-shortage equation may require an unrealistically low, or even negative, protection point to achieve the target. This condition occurred in 28 SKUs, mostly low-volume innovations.

3. **Excessive right-tail exposure:** For 21 SKUs, the computed SS exceeded one month of sales, implying an obsolescence risk.

These cases are outside the scope of the present study and may constitute a path for future research.

### 4.2.2 Protection-point calculation

In the literature, most authors work with CSL, although some examples of BFR can be found. Adenso-Díaz (1996) uses the classical normal-based fill-rate formulation and derives the expected shortage conditional on stockout occurrence. Silver et al. (1998, Chapter 7) provide a more general formulation for non-normal demand contexts that can be flexibly adapted to different distributions. Kanet et al. (2010) determine SS by using the standard normal probability density and cumulative distribution functions. Janssens and Ramaekers (2011) use the expected number of units short as a service measure when the lead-time demand distribution is only partially known. Jonsson and Mattsson (2019) warn of the method's limitations and simulate the entire system to then adjust an input parameter until the service target is reached.

Following Silver et al. (1998, Chapter 7), the expected shortage associated with a protection point under a continuous forecast-error distribution with density is:

$$E[(X-S)^{+}] = \int_{S}^{\infty} (x-S) f(x)\, dx$$

In this study, $X$ represents the CFEAP and $S$ is the protection buffer added to the forecast. For a target fill rate $\beta$ and replenishment lot size $Q$, the protection point is obtained from:

$$E[(X-S)^{+}] = (1-\beta_{used})Q$$

For each SKU, the forecast-error distribution is fitted parametrically:

$$X \sim D(\theta)$$

where $D$ denotes the selected distributional family and $\theta$ its estimated parameters. The required protection buffer $S^{*}$ is therefore defined as the numerical solution of:

$$\int_{S^*}^{\infty} (x - S^*) f(x; \theta)\, dx = (1 - \beta_{used})Q$$

The equation is solved numerically for each SKU using a bracketed root-finding procedure. The operating SS is then obtained from the resulting protection buffer. This formulation allows the protection point to reflect the fitted asymmetry, tail behaviour, and scale of the cumulative forecast-error distribution without imposing normality.

### 4.2.2 Beta Fill Rate objective and Cycle-End Exposure

Without SS, a normally distributed forecast error produces a stockout in half of the replenishment cycles. Under a BFR criterion, however, the expected stockout quantity will equal the standard deviation of demand during the lead time. As a result, the same service target can often be achieved with less SS. This advantage also makes the calculation more sensitive to persistent forecast errors: when errors occur in consecutive periods, shortages may become larger and last longer. A situation recorded by CSL as one failed cycle can therefore have a much greater effect on the achieved fill rate.

Lot size also affects the result. As the replenishment lot increases, the expected shortage represents a smaller proportion of the total quantity supplied during the cycle. The model may therefore require less SS to meet the same BFR target. This creates what we define as Cycle-End Exposure. Large lots produce long replenishment cycles and can make the average fill-rate result appear satisfactory even when protection becomes insufficient near the end of the cycle. The resulting stockout may then be operationally significant and, in some cases, last for more than one week, potentially breaching retailer service requirements and, in severe or recurrent cases, increasing the risk of the SKU being delisted.

A correction was therefore introduced to limit Cycle-End Exposure by constraining expected shortage to weekly demand $D_S$.

First, define the objective:

$$\beta_{obj} = 0.98$$

The expected stockout per cycle implied by the fill rate is:

$$ES = Q(1 - \beta)$$

$$Q(1 - \beta_{obj}) > DS$$

If this condition is not satisfied, the adjustment condition is defined as:

$$\beta_{lim} = 1 - \frac{DS}{Q}$$

$$\beta_{used} = \begin{cases} \beta_{lim}, & \text{if } Q(1-\beta_{obj}) > DS \\ \beta_{obj}, & \text{if } Q(1-\beta_{obj}) \leq DS \end{cases}$$

$$ES = Q(1 - \beta_{used})$$

## 4.3 Heteroscedasticity: FSEDM

FSEDM addresses the underexplored parametric modelling of forecast-dependent heteroscedasticity in safety-stock decisions. It preserves the SKU-specific distributional shape while allowing its scale to vary with the forecast level, so asymmetry and relative tail behaviour remain unchanged while weekly dispersion adapts to the estimated variance dependence. Visually, the method can be interpreted as a horizontal expansion or contraction of a common SKU-specific base distribution. The resulting conditional distribution is passed directly to the APFR expected-shortage calculation.

**Forecast-dependent heteroscedastic scale estimation**

For each SKU $i$, the historical forecast error is defined as:

$$E_{i,t} = VT_{i,t} - FC_{i,t}$$

where $VT_{i,t}$ is the observed sales and $FC_{i,t}$ is the forecast.

The proposed method assumes that the conditional variance of the forecast error is a linear function of the forecast:

$$\mathrm{Var}(E_{i,t} \mid FC_{i,t}) = \sigma_{i,t}^2$$

$$\sigma_{i,t}^2 = \omega_i + \beta_{\mathrm{HX},i} FC_{i,t}$$

where:

$$\omega_i \geq 0$$

$$\beta_{\mathrm{HX},i} \geq 0$$

The parameter $\omega_i$ represents the baseline component of error variability that is not driven by the forecast level. The parameter $\beta_{\text{HX},i}$ measures how strongly the error variability increases with the forecast.

In the implemented algorithm, these two parameters are estimated by regressing the squared forecast errors on the forecast values, with non-negativity constraints. For each SKU $i$, the estimation problem is:

$$(\widehat{\omega}_i, \hat{\beta}_{\text{HX},i}) = \arg \min_{\omega_i, \beta_{\text{HX},i}} \sum_{t=1}^{T_i} \left(E_{i,t}^2 - \omega_i - \beta_{\text{HX},i} FC_{i,t}\right)^2$$

subject to:

$$\omega_i \geq 0$$
$$\beta_{\text{HX},i} \geq 0$$

This is estimated using non-negative least squares. The non-negativity constraints are necessary because the model represents a conditional variance, which cannot be negative.

Once $\widehat{\omega}_i$ and $\hat{\beta}_{\text{HX},i}$ have been estimated, the weekly conditional scale is computed as:

$$s_{i,t} = \sqrt{\widehat{\omega}_i + \hat{\beta}_{\text{HX},i} FC_{i,t}}$$

The parameter $\hat{\beta}_{\text{HX},i}$ determines the sensitivity of the error-distribution width to the forecast level.

Then the method separates the forecast-error distribution into two components:

$$E_{i,t} = s_{i,t} Z_i$$

where $Z_i$ is the standardized base error distribution for SKU $i$, and $s_{i,t}$ is the forecast-dependent scale.

If the original BDFT-fitted distribution for the forecast error $E_i$ is described by:

$$E_i \sim F_E(\theta_{E,i})$$

with parameters:

$$\theta_{E,i} = (\text{shape parameters}, loc_{E,i}, scale_{E,i})$$

then the method defines a reference forecast:

$$FC_{\mathrm{ref},i} = \mathrm{median}(FC_{i,t})$$

The reference scale is:

$$s_{\mathrm{ref},i} = \sqrt{\widehat{\omega}_i + \hat{\beta}_{\mathrm{HX},i} FC_{\mathrm{ref},i}}$$

The standardized base distribution is then obtained by:

$$Z_i = \frac{E_i}{s_{\mathrm{ref},i}}$$

For location–scale distributions, this transformation implies:

$$loc_{Z,i} = \frac{loc_{E,i}}{s_{\mathrm{ref},i}}$$

$$scale_{Z,i} = \frac{scale_{E,i}}{s_{\mathrm{ref},i}}$$

The shape parameters are kept unchanged:

$$a_{Z,i} = a_{E,i}$$

$$b_{Z,i} = b_{E,i}$$

$$c_{Z,i} = c_{E,i}$$

For each week $t$, the conditional forecast-error distribution is obtained by rescaling the standardized base distribution:

$$E_{i,t} = s_{i,t} Z_i$$

Equivalently, if $f_{Z_i}(z)$ is the density of the standardized error distribution, the conditional density of the original forecast error is:

$$f_{E_i|FC}(e \mid FC_{i,t}) = \frac{1}{s_{i,t}} f_{Z_i}\left(\frac{e}{s_{i,t}}\right)$$

Therefore, BDFT continues to determine the distributional shape, while the conditional weekly scale $s_{i,t}$ adjusts its width according to the forecast. Visually, this is equivalent to stretching or compressing the same curve along the X-axis: as the forecast volume increases, $s_{i,t}$ increases and the 98% point moves farther from the centre, without altering the relative skewness or tail behaviour. Figure 2 illustrates this mechanism for three forecast levels of the same SKU, showing how the conditional error distribution widens while preserving the shape selected through BDFT.

**Figure 2. Forecast-scaled conditional error distributions under FSEDM**

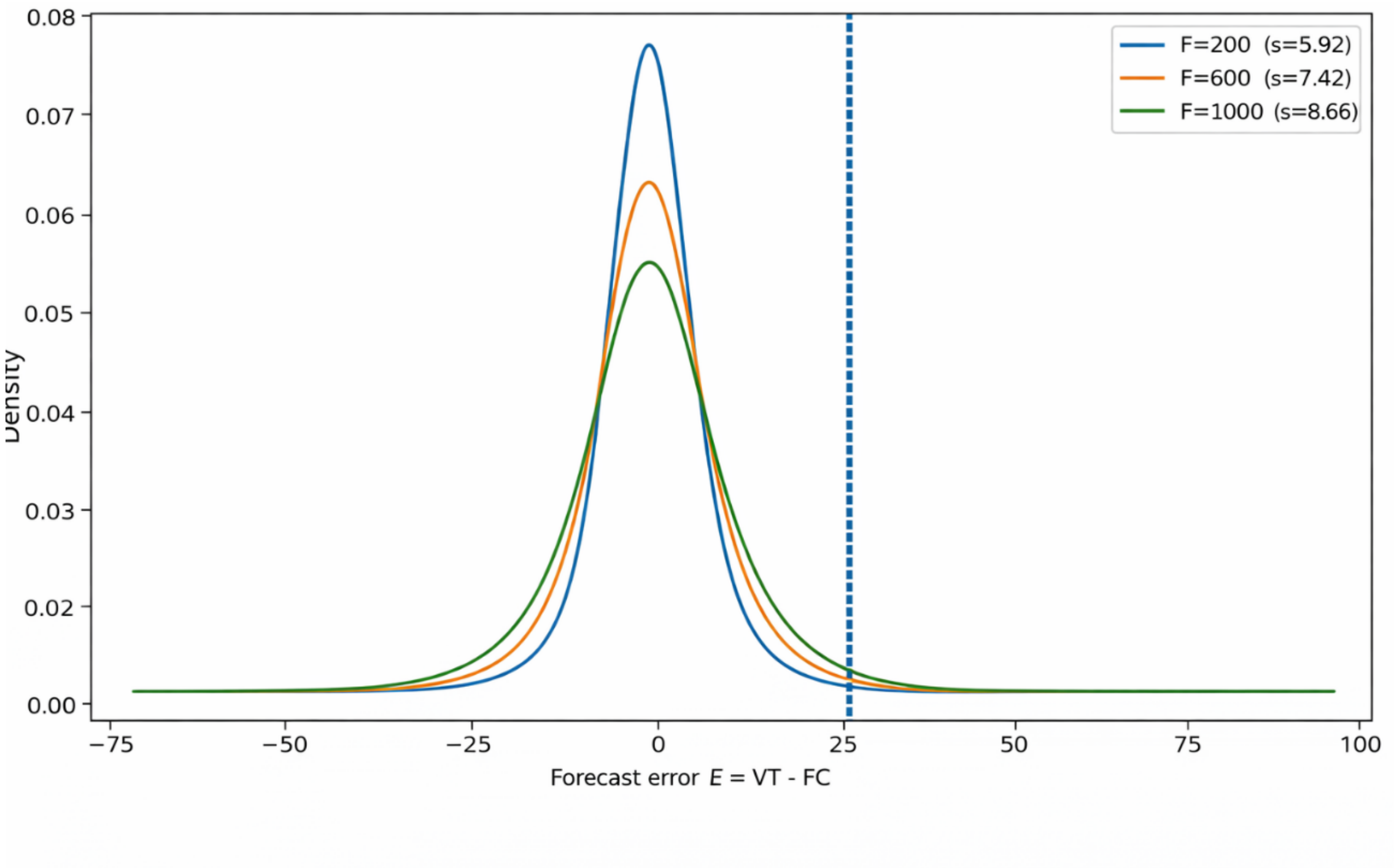

# 5. Simulation design

## 5.1 Simulation design and benchmark methods

**Benchmark method.** The study addresses a real industrial problem and therefore benchmarks the proposed approach against the company's current safety-stock method. As no single existing method jointly addresses the three identified deficiencies—forecast-error autocorrelation, conditional heteroscedasticity, and heavy-tailed behavior—the current company method provides the operational baseline against which the incremental effect of the integrated corrections can be assessed.

**Database and forecast-error construction.** Both methods are evaluated using the same weekly forecast-error database. Building on Fernández-Palacios et al. (2026), LOWDII is retained as a representativeness filter because it had already been validated on the same dataset. In the present study, it is applied to both methods so that the comparison isolates the contribution of the new methodological corrections. For the company's current method, LOWDII is applied directly to the raw weekly forecast errors; for the proposed method, it is applied after CFEAP transformation to represent cumulative uncertainty over the operational action period.

**Calibration and validation windows.** The period 2021–2024 is used for calibration because four complete calendar years provide a representative basis for estimating forecast-error uncertainty at high service levels (Fernández-Palacios et al., 2026), while 2025 is reserved for out-of-sample validation through discrete-event simulation as described in Section 5.2.

**Experimental unit.** The SKU is used as the experimental unit. The final simulation sample contains 153 SKUs.

**Performance metrics.** The methods are evaluated using four simulation outputs: average weekly inventory (AvST), beta fill rate (BFR), percentage of weeks without stockout (PctWeeks), and average stock out duration (AvOutDays). These metrics are defined in Section 5.3.

**Statistical comparison.** Since both methods are evaluated on the same SKUs, performance is compared using paired t-tests and paired Cohen's d for each metric. A MANOVA assesses the joint difference across all performance indicators. The incremental effect of FSEDM is further examined using group comparisons and HC3-robust regressions based on the intensity of heteroscedasticity. The rationale is detailed in Section 5.3 and the interpretation of these tests is detailed in Section 6.3.

## 5.2 Discrete-event simulation model

The discrete-event simulation (DES) evaluates the alternative safety-stock methods under the company's operational planning logic. All methods are evaluated under exactly the same demand sequence, forecast information, production constraints, lot sizes, capacity-control rule, event ordering, and shortage assumptions. The only systematic difference among simulation treatments is therefore the safety-stock policy generated by the alternative forecast-error calibration methods.

The DES represents finished-goods inventory at SKU level using weekly time buckets and reproduces the main operational constraints that interact with safety-stock decisions: a four-week production lead time, rolling forecast updates, minimum production lots, production increments, and shared production-capacity groups.

**Planning and replenishment logic**

The simulation uses a weekly state-transition model. For SKU $i$ and week $t$, inventory evolves according to:

$$ST_{i,t+1} = ST_{i,t} + PD_{i,t} - D_{i,t},$$

where $ST_{i,t}$ is opening inventory, $PD_{i,t}$ is production received during week $t$, and $D_{i,t}$ is realized demand. Since the model uses weekly buckets, each week constitutes one decision epoch, with production receipts processed before demand realization. No intra-week state transitions are represented; consequently, undershoot arising from within-week demand or production timing is outside the model scope. A future extension could evaluate a five-week protection interval to capture this effect (Campbell, 1995).

The simulation applies an order-up-to policy over the production lead-time horizon. Production decisions are taken using the forecast available at the corresponding planning origin. With a production lead time of $L = 4$ weeks, projected stock for a target week $t$ is computed as:

$$STPROY_{i,t} = ST_{i,t-L} + \sum_{k=0}^{L-1} P\,D_{i,t-L+k} - \sum_{k=0}^{L} F\,C_{i,t-L,t-L+k},$$

where $FC_{i,t-L,t-L+k}$ denotes the forecast for week $t - L + k$ available at planning origin $t - L$. Thus, projected stock incorporates opening inventory, production already fixed within the frozen horizon, and the complete forecast profile from the planning origin through the target week. This formulation corresponds directly to the implemented four-week projected-stock calculation.

The inventory gap for SKU $i$ in week $t$ is then defined as:

$$GAP_{i,t} = SSt_{i,t} - STPROY_{i,t},$$

where $SSt_{i,t}$ is the safety-stock target generated by the method being evaluated. Positive values of $GAP_{i,t}$ indicate that additional production is required to restore the projected inventory position to the desired protection level.

Production quantities respect both the SKU-specific minimum lot $Q_i$ and the production increment $INC_i$. When $GAP_{i,t} > 0$, the required production quantity is:

$$PD_{i,t}^{req} = Q_i + \left\lceil \frac{\max(Q_i, GAP_{i,t}) - Q_i}{INC_i} \right\rceil INC_i,$$

whereas no additional production is required when $GAP_{i,t} \leq 0$.

When several SKUs compete for the same production resource, the allocation sequence prioritizes the most critical projected inventory positions. SKUs with non-positive projected stock are processed first, in ascending order of $STPROY$, while the remaining SKUs are ranked by their uncovered safety-stock gap. Production is assigned only when the required batch fits within the available capacity of the corresponding production group, thereby preserving the interaction between safety stock, lot sizing, lead time, and capacity constraints documented for stochastic multi-item systems (Altendorfer, 2019).

**Figure 3. Weekly planning and inventory-transition logic used in the DES.**

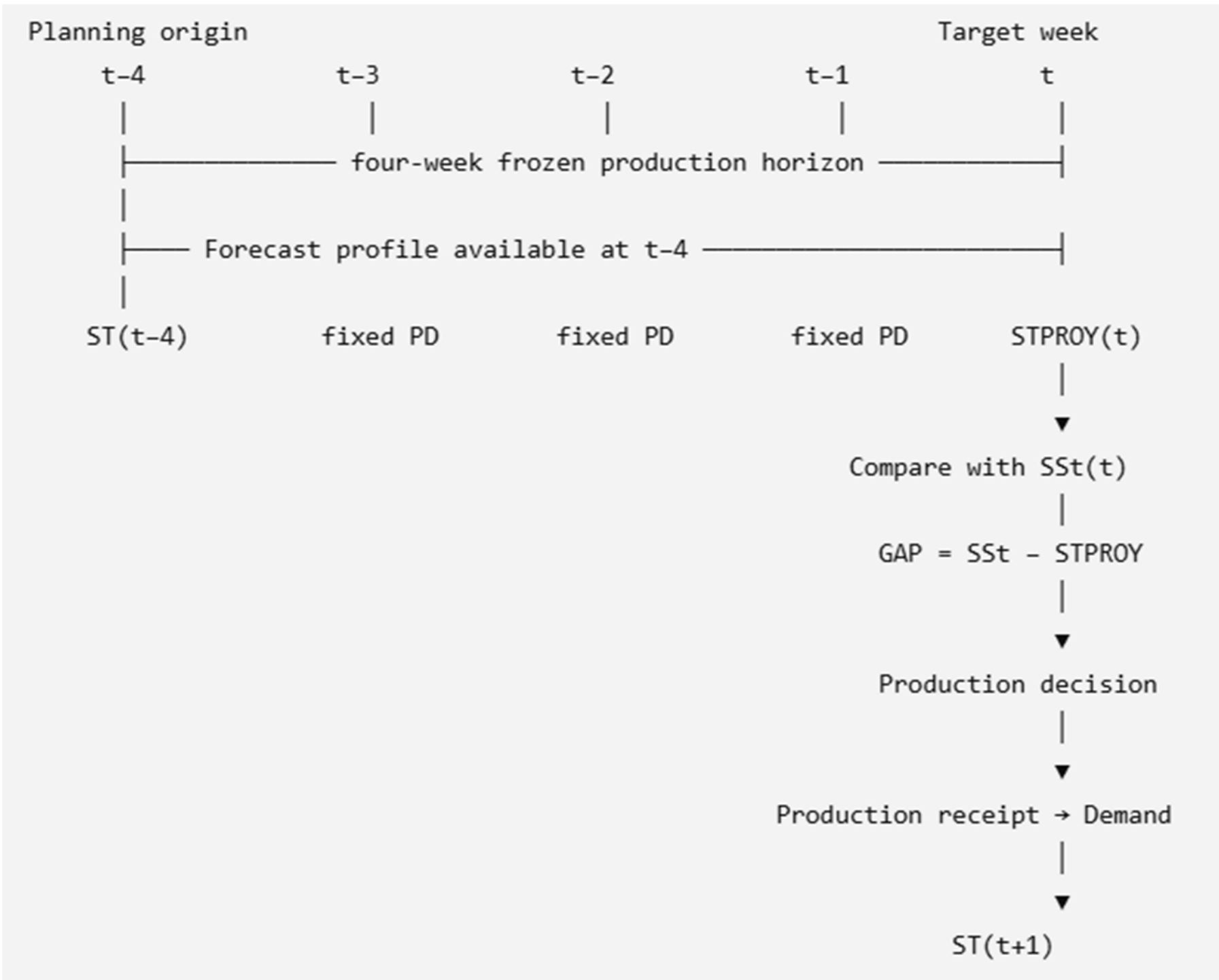

**Simulation horizon, warm-up, and repeated annual blocks**

The simulation starts in week 50 of 2024. The period from week 50 of 2024 to week 01 of 2025 is treated as a warm-up interval and is excluded from performance measurement. Production already fixed in the frozen horizon is preserved during this period so that the dynamic inventory state evolves consistently into the evaluation window.

The primary out-of-sample trajectory consists of the 52 weeks from week 02 of 2025 through week 01 of 2026. The model subsequently propagates the system state over repeated 52-week annual blocks. These blocks reuse the same 2025 exogenous demand, forecast, and safety-stock template, whereas endogenous inventory, production, backlog, and capacity carry-over states continue to evolve from one block to the next.

These repeated annual blocks do not constitute independent statistical replications and are not interpreted as additional out-of-sample observations. Their purpose is instead to reduce the dependence of low-rotation SKU performance on the particular initial position of the replenishment cycle. This is relevant because large production batches may cover several weeks of demand, so the effect of a demand disturbance may depend strongly on the inventory position inherited from previous production decisions.

**Scope and exclusions**

The simulation is restricted to finished-goods inventory planning. Component-level material constraints and upstream supply interactions are excluded to isolate the effects of forecast-error-based safety-stock calibration.

Only internally manufactured SKUs are included, representing 93.7% of total sales volume. Imported and co-packed products are excluded because their replenishment logic differs from the production system represented in the DES and because they account for only 2.4% of sales volume. Innovation items are also excluded, as their limited historical record does not provide sufficient observations for stable forecast-error calibration (Section 4.2.3); these items represent a further 3.9% of sales volume.

Although the company operates several production and logistics sites, inventory decisions are represented at a single location. This abstraction is consistent with recent simulation approaches, such as Disney et al. (2025), and avoids confounding the comparison of historical error treatments with multi-site allocation rules.

**Controlled-capacity design**

The simulation uses a controlled-capacity environment. For each production-resource group $g$ and week $t$, the baseline capacity envelope is defined as:

$$CapBase_{g,t} = \sum_{i \in g} D_{i,t} \, .$$

Thus, the baseline production envelope is normalized to the realized aggregate demand of the SKUs sharing that resource. This prevents structural capacity shortages from confounding the comparison among safety-stock methods. Consistently, the model excludes rationing behaviour associated with widespread supply crises. Shortages are instead assumed to be distributed randomly across a large customer base, allowing customer-level variability to be partially absorbed through risk pooling (Lee et al., 1997).

The actual capacity available to the weekly production-allocation routine additionally reflects changes in the aggregate safety-stock target and unused capacity carried forward from the preceding week:

$$Cap_{g,t} = \max\left[0, CapBase_{g,t} + \Delta SSt_{g,t} + \max(Carry_{g,t}^{in}, 0)\right].$$

As a result, realized demand is used only to construct the common ex-post capacity envelope, whereas replenishment decisions rely exclusively on forecasts available at the corresponding planning origin.

**Shortages and service assumptions**

A stockout occurs when the inventory position becomes negative after production receipts and realized demand have been processed for the week:

$$OUT_{i,t} = \max(-ST_{i,t}, 0).$$

The model assumes backordering rather than lost sales. Consequently, unmet demand remains embedded in the negative inventory position and must be recovered by subsequent production.

The DES should therefore be interpreted as a controlled comparative experiment rather than as a complete digital replica of the production system.

## 5.3 Performance metrics and statistical comparison tests

We measure the performance of our proposal through two primary indicators:

- **Average weekly inventory (Un)**: The objective is to minimize the average level of inventory. AvST is expressed in Un, a standardized inventory unit used consistently throughout the analysis.
- **Beta Fill Rate (%):** BFR is the service indicator used by the company under study. It is computed as the ratio between units served and units requested. In this study, it is evaluated both at aggregate level and at SKU level to determine whether the simulated inventory policy achieves the target service level of 98%.

Because BFR does not capture the frequency or duration of shortages, two secondary indicators are included: the percentage of weeks without stockouts and the average stockout duration. These metrics identify whether the target fill rate is achieved through frequent small shortages or fewer but more persistent disruptions.

- **Percentage of weeks without stockout (%):** PctWeeks is included to assess whether targeting BFR increases the frequency of stockout periods, even when the proportion of unserved units remains limited. This is particularly relevant for high-rotation SKUs, where frequent micro stockouts may coexist with a high fill rate.
- **Average stockout duration (days):** AvOutDays is included specifically to assess whether Cycle-End Exposure leads to prolonged shortages. Several short, low-volume stockouts are operationally different from a single extended period with no supply, particularly when it overlaps with a jointly funded commercial action.

For each indicator, we apply a paired t-test to assess whether the average difference between the proposed method and Current method is statistically different from zero. We also report paired Cohen's d to quantify the magnitude of the differences and distinguish statistical significance from practical relevance. In addition to the univariate analysis, we perform a MANOVA on the vector of SKU-level differences across the four dependent variables. This joint analysis evaluates whether the two methods differ overall when inventory, service, and stockout performance are considered simultaneously. Standard MANOVA statistics, including Wilks' Lambda, Pillai's Trace, Hotelling–Lawley Trace, and Roy's Greatest Root, are used to assess whether the joint vector of mean differences departs from zero.

To assess whether the incremental effect of FSEDM depends on the degree of heteroscedasticity, method differences are calculated at SKU level for both AvST and BFR. These differences are first compared between homoscedastic and heteroscedastic SKUs using a binary group model. Because heteroscedasticity is inherently continuous, additional regressions use the Breusch–Pagan LM statistic as a continuous measure of its intensity, avoiding the loss of information and statistical power caused by dichotomization. MANOVA is also applied to evaluate whether the joint stock–service effect varies with heteroscedasticity. Robust HC3 standard errors are used in the regression models to limit the influence of non-constant residual variance.

# 6. Results

## 6.1 Aggregate performance results

At the aggregate level, Table 1 reports portfolio-level results across the 153 SKUs included in the simulation, with average metrics weighted by sales volume where applicable. The proposed method achieves a BFR 0.68 percentage points higher than company's current method, although still below the 98% target. At the same time, it reduces average global SS by 8.1% and AvST by 9.9%, indicating a simultaneous improvement in inventory efficiency and aggregate fill-rate performance. Importantly, LOWDII alone reduces AvST by 7% relative to the company's baseline formulation, indicating that distributional-influence filtering can improve inventory efficiency while preserving the operational service target. The larger reduction obtained with the proposed method is consistent with an additional contribution from the remaining methodological components.

The 19.5% reduction in maximum stock is particularly relevant given the company's inflexible storage capacity, suggesting that the current seasonal method generates demand-related inventory peaks without commensurate service benefits.

Regarding the secondary metrics, average stockout duration decreases from 2.34 to 1.83 days, corresponding to a 21.8% reduction at the aggregate demand-weighted level. Although the SKU-level paired analysis reported below yields a small positive mean difference, this difference is not statistically significant. Hence, the results provide no statistically detectable evidence that the proposed method increases stockout duration or induces the hypothesized Cycle-End Exposure effect.

The percentage of weeks without stockout decreases by 0.8 percentage points, from 96.1% to 95.3%. This moderate decline is consistent with the BFR-based design of the proposed method, which controls the expected quantity of unmet demand rather than maximizing the probability of avoiding a stockout in every week. Given its limited magnitude and the absence of deterioration in BFR or stockout duration, the result does not indicate a material operational concern.

**Table 1. Aggregate level results comparing both methods**

| | Proposed Method | Current Company Method | Absolute difference | Relative difference (%) |
|---|---|---|---|---|
| **Average weekly SS (Un)** | 183,912 | 200,174 | -16,262 | -8.1% |
| **Average weekly inventory (Un)** | 259,334 | 287,704 | -28,370 | -9.9% |
| **Maximum weekly inventory (Un)** | 298,190 | 370,605 | -72,415 | -19.5% |
| **Beta fill rate (%)** | 97.36% | 96.69% | +0.68 pp | +0.7% |
| **Average stockout duration (days)** | 1.83 | 2.34 | -0.51 | -21.8% |
| **Weeks without stockout (%)** | 95.3% | 96.1% | -0.8 pp | -0.9% |

## 6.2 Statistical significance of performance differences

The SKU-level comparison provides clear evidence that the proposed method reduces inventory. AvST is significantly lower ($p = 0.0022$), with an estimated mean reduction of 185.4 units per SKU. The standardized effect is small to moderate ($d = -0.251$), indicating that the reduction is meaningful but modest relative to the substantial variability across SKUs. Operationally, this result confirms that the aggregate inventory reduction is not driven solely by a small number of high-volume products; it is also observable across the SKU portfolio.

For BFR, the proposed method records an average SKU-level difference of $-0.31$ percentage points relative to the current method. This difference is small ($d = -0.131$) and not statistically significant ($p = 0.1079$). Therefore, the analysis does not provide sufficient evidence that the inventory reduction generated by the proposed method causes a systematic deterioration in the principal service metric.

The SKU-level result differs from the aggregate weighted BFR, which is higher under the proposed method. This reversal arises because the two indicators apply different weighting structures. The paired analysis assigns equal weight to every SKU, whereas aggregate BFR is weighted by the number of units demanded, giving greater influence to high-volume products. Consequently, a slightly lower unweighted mean can coexist with a higher portfolio-level weighted result. This constitutes a form of Simpson's paradox, in which the direction of a comparison changes after observations are aggregated using different group weights (Simpson, 1951). The relevant conclusion is therefore not that both methods produce identical service, but that no statistically significant SKU-level BFR deterioration is detected, while the volume-weighted portfolio indicator favours the proposed method.

AvOutDays shows only a small and statistically non-significant change at SKU level. The estimated mean paired difference is 0.205 days per SKU, with a small effect size ($d = 0.094$) and no statistically significant difference ($p = 0.249$). Thus, there is no evidence that the inventory reduction achieved by the proposed method systematically increases the duration or operational severity of stockout episodes across SKUs. This result is consistent with the combined effect of the cumulative forecast-error treatment and the adjustment introduced to limit Cycle-End Exposure, although the present comparison does not isolate the contribution of each component separately.

PctWeeks, which measures the proportion of periods without a stockout, is 0.83 percentage points lower under the proposed method. The difference is statistically significant ($p = 0.0233$), but its standardized magnitude is small ($d = -0.185$). The proposed method therefore produces a slight reduction in the frequency with which complete stockout avoidance is achieved. However, this result has secondary operational importance because the study's primary service objective is BFR rather than Cycle Service Level. Considered together with the absence of a statistically significant change in BFR and AvOutDays at SKU level, the result suggests a small increase in the occurrence of stockout periods without evidence of a material increase in the proportion of unserved units or in shortage persistence.

The joint MANOVA confirms that the overall performance profiles of the two methods differ (Wilks' $\lambda = 0.918, p = 0.0135$). This multivariate result is consistent with the combined pattern observed in the individual metrics: a clear inventory reduction, a small decline in PctWeeks, and no statistically detectable deterioration in BFR or AvOutDays. It should not be interpreted as showing that every performance dimension differs independently, but rather that the complete stock–service vector changes when the metrics and their covariance are considered jointly.

The additional analysis of the Forecast-Scaled Error Distribution Method (FSEDM) shows no significant portfolio-wide average effect on stock or BFR. Likewise, dividing SKUs into binary homoscedastic and heteroscedastic groups does not reveal a significant differential effect. However, when heteroscedasticity is represented continuously by the Breusch–Pagan LM statistic, the benefit of FSEDM increases significantly with its intensity. Moving from the 25th to the 75th percentile of this statistic is associated with an expected BFR improvement of approximately 0.93 percentage points.

## 6.3 Interpretation: Statistical profiling and operational actionability

The proposed protection-point calculation does not reduce SS through a uniform scalar adjustment; instead, it recalibrates the buffer according to the empirical shape of the forecast-error distribution and the expected shortage condition associated with the fill-rate objective. The density-based calculation acts as a control mechanism: it protects against the relevant part of the right tail without imposing unnecessary protection where the fitted distribution does not support it. As the protection point is located more precisely, equivalent service performance can be achieved with significantly less stock, although with the expected trade-off of a slight increase in stockout frequency.

The results indicate that the proposed method improves inventory efficiency without transferring the reduction in stock into more severe service failures, as the simulation does not show a statistically significant deterioration in AvOutDays. This result is consistent with the intended role of the cumulative forecast-error treatment in mitigating consecutive stockouts associated with persistent forecast deviations, and with the Cycle-End Exposure adjustment in limiting the duration and severity of shortages occurring at the end of long replenishment cycles. However, the present comparison evaluates the combined method and does not isolate the individual contribution of either mechanism.

Finally, the results support the logic of the heteroscedastic extension, the method protects SKUs where forecast-dependent volatility is operationally relevant, while avoiding systematic overprotection in homoscedastic items. This selective treatment is consistent with the observed reduction in maximum stock, since inventory peaks are no longer amplified uniformly during high-forecast periods.

However, the aggregate results must be interpreted in light of an out-of-sample limitation. None of the methods reaches the 98.0% BFR target in 2025, mainly because two high-rotation SKUs experienced a distributional break relative to the 2021–2024 calibration period. Although they represent only two of the 153 SKUs included in the simulation, these items account for approximately 35% of the total stockout volume observed across the simulated portfolio. As shown in Figure 4, unforecast promotional activity generated large positive forecast errors for these SKUs, disproportionately reducing aggregate service. The shortfall therefore reflects an unmodelled change in the forecast-error process rather than a specific failure of the proposed calibration method.

**Figure 4. Distributional shift between calibration and validation forecast errors in two high-rotation SKUs**

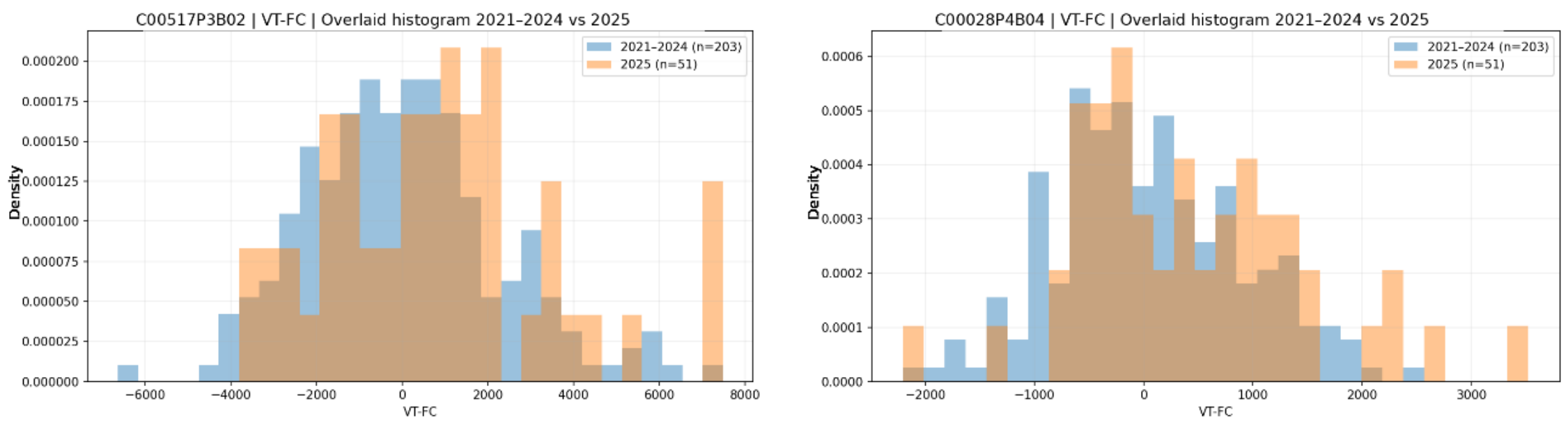

## 6.4 Managerial implications

The principal methodological contribution is an evidence-triggered architecture: solutions are not imposed uniformly across the portfolio but are activated, parameterized, or constrained according to the statistical profile and operational characteristics of each SKU. As a result, the proposed method avoids homogeneous treatment across SKUs and instead applies differentiated responses to statistically and operationally distinct cases.

Second, the proposed method addresses the main statistical limitations observed in the data. Its traceable architecture may improve managerial understanding and help mitigate algorithm aversion, although this effect is not directly evaluated.

Third, each adjustment has a clear operational interpretation, allowing managerial intervention to shift from subjective output correction to reviewing the statistical causes behind those outputs.

Fourth, the method is compatible with the company's existing ERP/MRP environment. This is important because the main barrier to adoption is often not statistical sophistication itself, but the cost and organizational disruption associated with implementation. The output can be expressed as an updated safety-stock parameter that can be loaded into the current planning system.

## 6.5 Limitations and future research directions

In practice, the company reached its 98% target service level not by adding incremental stock to the current standard, but through a specific operational risk-pooling structure, whereby shortages in some SKUs could be partially absorbed by available inventory elsewhere in the system. This effect is not captured in the present simulation, which isolates the finished-goods safety-stock decision at SKU level. Future work could incorporate this pooling structure to represent portfolio-level risk absorption, potentially improving stock efficiency while preserving the aggregate service objective.

Secondly, future research could compare the frequency and magnitude of manual planner overrides before and after implementation to assess whether the method's traceable architecture improves managerial trust and reduces algorithm aversion.

Finally, future research could examine alternative approaches to safety-stock dimensioning for innovations and other SKUs with insufficient or unstable historical records. In the present study, several SKUs could not be incorporated into the final simulation because the method could not identify a reliable distribution due to limited data, pronounced negative forecast-error bias, or excessively heavy-tailed error profiles. Most of these cases corresponded to low-volume innovations, for which the proposed method may require a different logic from that applied to mature SKUs. Future research could therefore examine specific procedures for short-history products, combining limited internal records with analogue-based estimation, lifecycle-stage information, or more conservative obsolescence-aware inventory rules.

# 7. Conclusions

This study develops and validates an adaptive safety-stock calibration method for environments characterized by temporally dependent, non-normal, biased, and forecast-dependent heteroscedastic errors, together with lot-size-related exposure under a BFR objective. Building on established statistical principles, the method contributes a coordinated sequence of conditionally activated modules that share a common error representation and translate SKU-specific statistical features into a single operational protection-point calculation.

The simulation results support the operational value of the proposed approach. Compared with the company's current seasonal safety-stock method, the proposed method reduces average weekly inventory by 9.9%, without causing statistically detectable deterioration in BFR. This suggests that the method recalibrates the buffer to protect the right tail only where supported by the fitted distribution. This efficiency gain comes with the expected trade-off of a slight increase in stockout frequency.

The cumulative forecast-error treatment, together with the adjustment introduced to limit Cycle-End Exposure, is consistent with the observed inventory reduction without an increase in stockout persistence, as measured by average stockout duration.

The method protects SKUs where forecast-dependent volatility is operationally relevant, while avoiding systematic overprotection in homoscedastic items. This produces one of the most notable optimisation effects of the method: a reduction of 19.5% in maximum inventory and greater operational stability.

Overall, the study shows that greater statistical realism can be combined with SKU-level adaptability, managerial interpretability, and operational implementability in safety-stock calibration, while reducing inventory investment and unnecessary inventory peaks without materially compromising service performance.

***Declaration of generative AI and AI-assisted technologies in the manuscript preparation process:*** *During the preparation of this work, the authors used Undermind to support the bibliographic search process and identify potentially relevant references. ChatGPT was used to assist in writing Python code for the simulation. All AI-assisted outputs were critically reviewed, verified, and edited by the authors. The authors take full responsibility for the accuracy, integrity, and final content of the manuscript.*